\documentclass[letterpaper, 10 pt, conference]{ieeeconf}  

\IEEEoverridecommandlockouts                              

\usepackage{times} 
\usepackage{amsmath} 
\usepackage{amssymb}  
\usepackage{amsfonts}
\usepackage[colorinlistoftodos,bordercolor=orange,backgroundcolor=orange!20,linecolor=orange,textsize=scriptsize]{todonotes}
\usepackage{hyperref}
\hypersetup{colorlinks, urlcolor = blue, linkcolor = blue, citecolor = red}
\usepackage{tikz}
\usetikzlibrary{shapes,arrows,decorations.markings,positioning}

\newtheorem{theorem}{Theorem}

\newtheorem{lemma}[theorem]{Lemma}
\newtheorem{corollary}[theorem]{Corollary}
\newtheorem{definition}[theorem]{Definition}
\newtheorem{remark}[theorem]{Remark}

\newtheorem{assumption}{Assumption}

\newtheorem{algo}[theorem]{Algorithm}

\newcommand{\R}{\mathbb{R}}

\newcommand{\E}{\mathcal{E}}
\newcommand{\Rmax}{\mathbb{R}_\text{max}}
\newcommand{\Emax}{\mathcal{E}_\text{max}}
\newcommand{\Rewards}{\mathcal{R}}
\newcommand{\C}{\mathcal{C}}
\newcommand{\G}{\mathcal{G}}
\newcommand{\F}{\mathcal{F}}
\newcommand{\scalar}[2]{\langle #1,#2 \rangle}
\newcommand{\MIN}{\textsf{MIN}}
\newcommand{\MAX}{\textsf{MAX}}
\newcommand{\Cpt}[1]{{S \setminus #1}}
\newcommand{\clo}[1]{{\bar{#1}}}

\DeclareMathOperator{\argmin}{arg\,min}
\DeclareMathOperator{\argmax}{arg\,max}

\newcommand{\unit}{\mathbf{1}}
\newcommand{\indic}{\mathbf{1}}
\newcommand{\sE}{\mathbb{E}}       
\newcommand{\ki}{k}
\newcommand{\kj}{\ell}
\newcommand{\set}[2]{\{#1\mid\,#2\}}

\renewcommand{\theenumi}{\roman{enumi}}

\title{\LARGE \bf
Generic uniqueness of the bias vector of mean payoff zero-sum games
}

\author{Marianne Akian, St\'ephane Gaubert and Antoine Hochart
\thanks{M. Akian, S. Gaubert and A. Hochart are with INRIA Saclay-Ile-de-France and CMAP, Ecole Polytechnique, CNRS, CMAP, Route de Saclay, 91128 Palaiseau Cedex, France.}
\thanks{The last author is supported by a PhD fellowship of Fondation Math\'ematique Jacques Hadamard.}
\thanks{{\tt\small marianne.akian@inria.fr} (M. Akian)}
\thanks{{\tt\small stephane.gaubert@inria.fr} (S. Gaubert)}%
\thanks{{\tt\small hochart@cmap.polytechnique.fr} (A. Hochart)}%
}

\begin{document}

\maketitle
\thispagestyle{empty}
\pagestyle{empty}

\begin{abstract}
Zero-sum mean payoff games can be studied by means of a nonlinear spectral problem.
When the state space is finite, the latter consists in finding an eigenpair $(u,\lambda)$ solution of $T(u)=\lambda \unit + u$ where $T:\R^n \to \R^n$ is the Shapley (dynamic programming) operator, $\lambda$ is a scalar, $\unit$ is the unit vector, and $u \in \R^n$.
The scalar $\lambda$ yields the mean payoff per time unit, and the vector $u$, called the {\em bias}, allows one to determine optimal stationary strategies. 
The existence of the eigenpair $(u,\lambda)$ is generally related to ergodicity conditions. 
A basic issue is to understand for which classes of games the bias vector is unique (up to an additive constant).
In this paper, we consider perfect information zero-sum stochastic games with finite state and action spaces, thinking of the transition payments as variable parameters, transition probabilities being fixed.
We identify structural conditions on the support of the transition probabilities which guarantee that the spectral problem is solvable for all values of the transition payments.
Then, we show that the bias vector, thought of as a function of the transition payments, is generically unique (up to an additive constant).
The proof uses techniques of max-plus (tropical) algebra and nonlinear Perron-Frobenius theory.
\end{abstract}


\section{Introduction}

\subsection{The mean payoff of a repeated game}

Zero-sum repeated games describe long term interactions between two agents (called players) with opposite interests.
We consider here perfect information games, in which each of the players chooses alternatively an action, being informed of the current state of the game and of the previous actions of the other player.
These choices determine an instantaneous payment, as well as the next state, either by a deterministic or a stochastic process.
We refer the reader to~\cite{NS03} for background on repeated games and stochastic games.

A zero-sum repeated stochastic game with perfect information is composed of
\\ --  a state space $S$.
\\ --  an action space $A_i\subset A$ for player \MIN, depending on $i\in S$, and included in a given set $A$;
at each stage, player \MIN\ chooses an action $a\in A_i$ knowing the current state $i$ and the past states and actions of both players.
\\ --  an action space $B_{i,a}\subset B$ for player \MAX, depending on $i\in S$ and on $a\in A_i$, and included in a given set $B$;
at each stage, player \MAX\ chooses an action $b\in B_{i,a}$, knowing the current state $i$, the last action $a \in A_i$ chosen by player \MIN\  and the past states and actions of both players.
\\ --  a transition payment $r_{i}^{a b}\in \R$ paid by player \MIN\ to player \MAX\ at each stage, given the current state $i\in S$ and the last actions of each of the players  $a\in A_i$ and $b\in B_{i,a}$.
\\ --  a transition probability $P_i^{a b} \in \Delta(S)$, depending on the same data, where $\Delta(S)$ is the set of probabilities over $S$;
when in state $i\in S$, if the actions $a,b$ have been played, the next state is chosen according to the probability $P_{i}^{ab}$.

The game is said to be {\em deterministic} if each transition probability $P_i^{a b}$ has only one state in its support, meaning that for every choice of $i\in S$, $a\in A_i$, and $b\in B_{i,a}$, we have $P_{i}(\{j\})=1$ for precisely one state $j\in S$.
Given the above parameters, and a strategy of each player (that is a selection rule for the sequence of actions at each step of the game), the payoff of the game is an additive function of the transition payments, that player \MIN\ intends to minimize and player \MAX\ intends to maximize.
In particular, the game with a fixed finite horizon $k$ consists of $k$ successive alternated moves of the players \MIN\ and \MAX, with a payoff
\begin{equation}
  \label{payofffinite}
  J^k_i\,=\, \sE \left[ \, \sum_{\ell = 0}^{k-1} r_{\xi_\ell}^{\alpha_\ell \beta_\ell} \,\right], 
\end{equation}
where $\sE$ denotes the expectation for the probability law of the process $(\xi_\ell,\alpha_\ell,\beta_\ell)_{\ell\geq 0}$ of state and actions of players \MIN\ and \MAX\ respectively, determined by the above transition probabilities and the strategies, given the initial state $i$.

Throughout the paper, we make the following assumptions:
\begin{assumption}[finite state and action spaces]
  \label{StandardAssumptions}
  \leavevmode
  \begin{enumerate}
  \item The state space $S$ is finite, say $S=\{1,\dots,n\}$.
  \item All the action spaces $A_i$ and $B_{i,a}$ are nonempty finite sets.
  \end{enumerate}
\end{assumption}

The analysis of the above game involves the {\em Shapley} operator $T:\R^n \to \R^n$, which is such that, for all $x=(x_i)_{i\in S}$,
\begin{equation}
  \label{eq:ShapleyOperator}
  [T(x)]_i = \min_{a \in A_i} \max_{b \in B_{i,a}} \big( r_i^{a b} + P_i^{a b} x \big),
\end{equation}
for all states $i\in S$. 
Note that Assumption~\ref{StandardAssumptions} implies that the $\min$ and $\max$ in~\eqref{eq:ShapleyOperator} are always attained.
Also, the elements $P$ of $\Delta(S)$ are seen as row vectors, $P=(P_j)_{j \in S}$, so that $Px$ means $\sum_{j=1}^n P_j x_j$.

Given an initial state $i\in S$, the game in horizon $k$ is known to have a value, denoted by $v^k_i\in \R$, which is equal to the minimum over all strategies of \MIN\ of the maximum over all strategies of \MAX\ of the payoff $J^k_i$ of~\eqref{payofffinite}.
In particular, the value vector $v^k=(v^k_i)_{1\leq i\leq n}\in\R^n$ satisfies the dynamic programming recursion
\begin{equation}
  \label{eq:RecursiveValue}
  v^0 = 0 \enspace \text{and} \enspace v^{k+1} = T(v^k) .
\end{equation}

Note that, since all the action spaces are finite, the above game can be transformed into a game with same behavior and operator~\eqref{eq:ShapleyOperator}, but with action spaces independent of the state and the selected action (by duplicating some transition payments and probabilities).
However, this would change the space of transition payments hence most of our results.

We will study here the asymptotics of the value vector as the horizon $k$ of the game tends to infinity.
In particular, we are interested in the following quantity, called the {\em mean payoff} vector (per time unit):
\begin{align}
  \label{e-def-chi}
  \chi(T) := \lim_{k \to \infty} v^k/k = \lim_{k\to \infty} T^k(0)/k \enspace .
\end{align}

The mean payoff $\chi(T)$ does exist when the action spaces are finite.
Indeed, Kohlberg~\cite{Koh80} has shown that $T$ has then an invariant half-line, that is, there exist two vectors $x, \nu \in \R^n$ such that $T(x+t\nu) = x+(t+1)\nu$ for every scalar $t$ large enough.
It can be checked that $\chi(T) = \nu$. 
Note that the existence of the mean payoff vector is also guaranteed more generally when $T$ is semi-algebraic, see Neyman~\cite{Ney03}.
Other conditions, in terms of the operator $(\alpha, x) \mapsto \alpha T(\frac{1-\alpha}{\alpha}x), \enspace \alpha > 0$, have been given by Rosenberg and Sorin~\cite{RS01}. 

It is convenient to note that the limit $\lim_{k\to\infty} T^k(x)/k$ is independent of the choice of $x\in \R^n$.
This follows readily from the elementary fact that $T$ is sup-norm nonexpansive.

\subsection{Nonlinear spectral problem}

The existence of the mean payoff vector is guaranteed if we can find a vector $u\in\R^n$ and a scalar $\lambda\in \R$ solution of the following nonlinear spectral problem,
\begin{equation}
  \label{eq:NonlinearSpectralProblem}
  T(u) = \lambda \unit + u,
\end{equation}
where $\unit$ is the unit vector of $\R^n$. 

Indeed, in such a situation, we have $T^k(u) = k\lambda +u$, for all $k\geq 0$, and so, $\chi(T) = \lambda \unit$, meaning that the mean payoff starting from state $i$, $[\chi(T)]_i$, is independent of the initial state $i$, and equal to $\lambda$.
The scalar $\lambda$ will be called the {\em eigenvalue} of $T$.
It follows that the eigenvalue, if it exists, is necessarily unique (for $\chi(T)$ is uniquely defined).
The vector $u$, called {\em bias vector} or {\em eigenvector}, gives optimal stationary strategies. 

Since any Shapley operator~\eqref{eq:ShapleyOperator} is additively homogeneous (that is, commutes with the addition of a constant), a bias vector is always defined up to an additive constant and it is interesting to understand when such a vector is unique (up to an additive constant). 

For one player problems, i.e., for discrete optimal control, the nonlinear eigenproblem~\eqref{eq:NonlinearSpectralProblem} (also known as the `average case optimality equation') has been much studied, either in the deterministic or in the stochastic case (Markov decision problems).
Then, the representation of bias vectors and their relation with optimal strategies is well understood.

For the deterministic case, the analysis relies on max-plus spectral theory, which goes back to the work of Romanovski~\cite{Rom67}, Gondran and Minoux~\cite{GM77} and Cuninghame-Green~\cite{CG79}.
In~\cite{MS92, BCOQ92, Bap98} can be found more background on max-plus spectral theory.
Kontorer and Yakovenko~\cite{KY92}, and Kolokoltsov and Maslov~\cite{KM97}, deal specially with infinite horizon optimization and mean payoff problems.
The set of bias vectors has the structure of a max-plus (tropical) cone, i.e., it is invariant by max-plus linear combinations, and it has a unique minimal generating family consisting of certain `extreme' generators, which can be identified by looking at the support of the maximizing measures in the linear programming formulation of the optimal control problem, or at the ``recurrence points'' of infinite optimal trajectories.
We refer the reader to~\cite{ABG07} and to the references therein for more information on max-plus spectral theory.
A combinatorial interpretation of some of these results, in terms of polyhedral fans, has been recently given by Sturmfels and Tran~\cite{ST13}.
The eigenproblem has also been studied for an infinite dimensional state space in the context of infinite dimensional max-plus spectral theory, see Akian, Gaubert and Walsh~\cite{AGW09}, and also in the setting of weak KAM theory, for which we refer the reader to the book of Fathi~\cite{Fat08}.
In the stochastic case, the structure of the set of bias vectors is still known, at least when the state space is finite, see Akian and Gaubert~\cite{AG03}. 

In the two player case, the structure of the set of bias vectors $u\in \R^n$ is less well known, although the description of this set remains a fundamental issue.
In particular, the uniqueness of the bias vector (up to an additive constant) is an important matter for algorithmic purposes.
Indeed, the nonuniqueness of the bias typically leads to numerical unstabilities or degeneracies.
In particular, the standard Hoffman and Karp policy iteration algorithm~\cite{HK66} may fail to converge in situations in which the bias vector is not unique.
Some refinements of the Hoffman and Karp scheme have been proposed by Cochet-Terrasson and Gaubert~\cite{cochet-cras}, Akian, Cochet-Terrasson, Detournay and Gaubert~\cite{ACTDG12}, and Bourque and Raghavan~\cite{bourque}, allowing one to circumvent such degeneracies at the price of a complexification of the algorithm (handling the non-uniqueness of the bias).
Hence, it is of interest to understand when such technicalities can be avoided.

Moreover, generally, the {\em existence} of the bias vector is controlled by some ergodicity conditions.
A general result in~\cite{GG04} relates the existence of the bias vector of a Shapley operator with the uniqueness of the bias vector of an auxiliary (simpler) Shapley operator, called recession function.
This provides further motivation to study the uniqueness problem for the bias vector.

\subsection{Main results}

We deal with two problems concerning zero-sum {\em parametric} stochastic games, that is, games whose parameters (transition payments and probabilities) vary.

We first look for conditions involving the transition probabilities which guarantee that the nonlinear eigenproblem is solvable {\em for all values of the transition payments}.
Our first main result, Theorem~\ref{thm:Support}, shows that this property is structural, meaning that it only depends on the support of the transition probabilities (not on their values).

Then, we address the question of the uniqueness of the bias vector, restricting our attention to games for which the previous structural property holds.
Our second main result, Theorem~\ref{thm:GenericUniqueness}, shows that the bias vector is generically unique (up to an additive constant).
More precisely, we show that the set of transition payments for which the bias vector is not unique belongs to a polyhedral fan the cells of which have codimension one at least.
A first ingredient in the proof is the max-plus spectral theorem, in the deterministic case, and its extension to the stochastic case~\cite{AG03}.
A second ingredient is a general result, showing that the set of fixed points of a nonexpansive self-map of $\R^n$ is a retract of $\R^n$, see Theorem~\ref{thm:NonexpansiveRetract}.
This allows us to infer the uniqueness of the bias vector of a Shapley operator from the uniqueness of the bias vector of the reduced Shapley operators obtained by fixing the strategy of one player.

Finally, as a simple application of our results, we show that the standard Hoffman and Karp policy iteration for two player games does terminate for generic well posed instances.

Let us note that although the finiteness of action spaces is assumed throughout the paper, the results in Section~\ref{sec:StructuralEigenvalue} remain true with compact action spaces (and some additional assumptions on the transition payments and probabilities, see~\cite{AGH14}).
In particular, they are true for repeated games with imperfect information (where players select an action simultaneously).
However, the techniques used to prove the results in Section~\ref{sec:GenericUniqueness} are strongly related with the finiteness of action spaces.
Therefore, they cannot go through Assumption~\ref{StandardAssumptions} and the perfect information assumption.

The paper is organized as follows.
In Section~\ref{sec:StructuralEigenvalue}, we show the structural existence result for the bias vector.
This partly relies on recent results established by the authors in a companion work~\cite{AGH14}, that we first recall.
The generic uniqueness of the bias vector is established in Section~\ref{sec:GenericUniqueness}, using results of max-plus spectral theory and nonlinear spectral theory.
The application to policy iteration is presented in Section~\ref{sec:Application}.

\section{Structural properties for the existence of an eigenvalue}
\label{sec:StructuralEigenvalue}

\subsection{Necessary and sufficient condition of existence of an eigenvalue}

Recall that, $\R^n$ being endowed with its usual partial order, the map $T:\R^n \to \R^n$ is said to be {\em monotone} and {\em additively homogeneous} if it satisfies the following properties, respectively,
\begin{itemize}
\item $x \leq y \enspace \Rightarrow \enspace T(x) \leq T(y), \enspace x,y \in \R^n$,
\item $T(x + \lambda \unit) = T(x) + \lambda \unit, \enspace x \in \R^n, \lambda \in \R$.
\end{itemize}
In~\cite{GG04}, Gaubert and Gunawardena have given a sufficient condition for a monotone and additively homogeneous operator $T$ on $\R^n$ to have an eigenvalue.
This condition involves the {\em recession function} of $T$, denoted by $\hat T$ and defined on $\R^n$ by
\[
  \hat T(x) := \lim_{\alpha \to +\infty} \frac{T(\alpha x)}{\alpha}.
\]
For a general operator, the recession function may not exist. However, if $T$ is a Shapley operator~\eqref{eq:ShapleyOperator} satisfying Assumption~\ref{StandardAssumptions} (finite state and action spaces), then the recession function is well defined everywhere and given by
\begin{equation}
  \label{eq:PaymentFreeOperator}
  [\hat T(x)]_i = \min_{a \in A_i} \max_{b \in B_{i,a}} P_i^{a b} x, \enspace \forall i \in S.
\end{equation}

\begin{remark}
  If the recession function $\hat T$ of $T$ exists, then it is also monotone and additively homogeneous.
  Furthermore, it is positively homogeneous, meaning that
  \begin{itemize}
  \item $\hat T(\lambda x) = \lambda \hat T(x), \enspace x \in \R^n, \lambda >0$.
  \end{itemize}
  As a consequence, any vector proportional to $\unit$ is a fixed point of $\hat T$.
\end{remark}

\begin{theorem}[\cite{GG04}]
  \label{thm:RecessionFunction}
  Let $T:\R^n \to \R^n$ be a monotone and additively homogeneous operator.
  Suppose its recession function is well defined on $\R^n$ and has only trivial fixed points (that is, fixed point proportional to $\unit$).
  Then $T$ has an eigenvalue, meaning that Equation~\eqref{eq:NonlinearSpectralProblem} is solvable for $T$.
\end{theorem}

Conversely, if $F$ is a monotone, additively and positively homogeneous operator on $\R^n$ and if $\nu$ is one of its fixed points, then the operator $T: x \mapsto F(x)+\nu$ is monotone, additively homogeneous and satisfies $\hat T = F$ and $\chi(T) = \nu$.
In particular, if $\nu$ is nontrivial, then $T$ cannot have any eigenvalue.
This leads to the following theorem.
\begin{theorem}
  \label{thm:StructuralEigenvalue}
  Let $T$ be a Shapley operator~\eqref{eq:ShapleyOperator} with finite state and action spaces (Assumption~\ref{StandardAssumptions}).
  Then the following properties are equivalent:
  \begin{enumerate}
  \item $T$ has an eigenvalue for all values of the transition payments $r_i^{a b}$,
  \item the recession function of $T$~\eqref{eq:PaymentFreeOperator} has only trivial fixed points.
  \end{enumerate}
\end{theorem}

\begin{remark}
  The special case of Theorem~\ref{thm:StructuralEigenvalue} restrained to min-max functions (which corresponds to Shapley operators of deterministic games) has been addressed by Zhao, Zheng and Zhu~\cite{ZZZ01}.
  In this case, the condition on the recession function can be reduced to a fixed point problem for a Boolean function (defined on the set $\{0,1\}^n$ of Boolean vectors), hence a purely combinatorial problem.
\end{remark}

\subsection{The condition of existence of an eigenvalue is structural}

\subsubsection{Statement of the result}
In the following section, we shall consider operators~\eqref{eq:ShapleyOperator} that have an eigenvalue for all values of the transition payments.
According to Theorem~\ref{thm:StructuralEigenvalue}, these are characterized by the property that the operator given by~\eqref{eq:PaymentFreeOperator} has only trivial fixed points.

Fixed point problems related with such operators~\eqref{eq:PaymentFreeOperator}, which we call {\em payment-free} Shapley operators, have been studied in our companion work~\cite{AGH14}.
Our first main result is the following.
\begin{theorem}
  \label{thm:Support}
  Whether a given payment-free Shapley operator~\eqref{eq:PaymentFreeOperator} has only trivial fixed points depends only on the support of the transition probabilities.
\end{theorem}

Recall that the support of a probability $P \in \Delta(S)$ is the subset defined as $\{j \in S \mid P_j \neq 0 \}$.
In other words, changing the probabilities $P_{i}^{ab}$, while leaving invariant the set of $(i,j,a,b)$ such that $P_{ij}^{ab}\neq 0$ does not change the existence of a non-trivial fixed point of the operator $\hat{T}$ in~\eqref{eq:PaymentFreeOperator}.

We next recall some constructions and results that we use to establish Theorem~\ref{thm:Support}.

\subsubsection{Galois connection between invariant faces of $[0,1]^n$}
\label{sec:GaloisConnection}
Let $F:\R^n \to \R^n$ be a payment-free Shapley operator.
If $u$ is a fixed point of $F$ then, denoting $I:=\argmin u$ and $J:=\argmax u$, we can show that
\begin{equation}
  \label{eq:InvariantFaces}
  F(\indic_\Cpt{I}) \leq \indic_\Cpt{I} \quad \text{and} \quad \indic_J \leq F(\indic_J),
\end{equation}
where $\indic_K$, with $K \subset S$, is the vector with entries $1$ on $K$ and $0$ on $\Cpt{K}$.

Let $\F^-$ (resp.\ $\F^+$) be the family of subsets of $S$ verifying the first (resp.\ the second) inequality in~\eqref{eq:InvariantFaces}:
\begin{align*}
  \F^- &:= \big\{ I \subset S \mid F(\indic_\Cpt{I}) \leq \indic_\Cpt{I} \big\},\\
  \F^+ &:= \big\{ J \subset S \mid \indic_J \leq F(\indic_J) \big\}.
\end{align*}
The elements of $\F^-$ and $\F^+$ can be visualized as faces of the hypercube $[0,1]^n$ that are invariant by $F$.

The families $\F^-$ and $\F^+$ are nonempty lattices of subsets (they both contain $\emptyset$ and $S$).
We define a Galois connection (see~\cite{Bir95} for the definition) between them, denoted by $(\Phi, \Phi^\star)$, in the following way: to $I \in \F^-$ we associate the greatest element $J \in \F^+$ which has an empty intersection with $I$, and vice versa.
Formally we have, for $I \in \F^-$ and $J \in \F^+$,
\begin{equation}
  \label{eq:GaloisConnection}
  \Phi(I) := \bigcup_{J \in \F^+: I \cap J = \emptyset} J \enspace \text{and} \enspace \Phi^\star(J) := \bigcup_{I \in \F^-: I \cap J = \emptyset} I.
\end{equation}
If $I \in \F^-$, we denote by $\clo{I} := \Phi^\star \circ \Phi(I)$ its closure by the Galois connection (likewise, $\clo{J}$ is the closure of $J \in \F^+$).

The following results explain why this Galois connection is useful for the problem of existence of nontrivial fixed points.
\begin{theorem}[\cite{AGH14}]
  \label{thm:NontrivialFP}
  Let $F$ be a payment-free Shapley operator and $I \in \F^-$.
  \begin{enumerate}
  \item If $\Phi(I)=\emptyset$, then $F$ has no nontrivial fixed point $u$ such that $\argmin u = I$.
  \item If $I \neq \emptyset$ and $I = \clo{I}$, then $F$ has a fixed point $u$ satisfying $\argmin u = I$.
  \end{enumerate}
\end{theorem}

\begin{corollary}
  \label{coro:GaloisConnectionFixedPoint}
  A payment-free Shapley operator has a nontrivial fixed point if and only if there is a nontrivial subset $I \in \F^-$ closed for the Galois connection $(\Phi, \Phi^\star)$, that is a set $I \in \F^- \setminus \{\emptyset, S\}$ such that $I = \clo{I}$.
\end{corollary}

\subsubsection{Boolean abstractions of Shapley operators}
Let $F$ be a payment-free Shapley operator.
We call {\em upper} and {\em lower Boolean abstractions} of $F$ the operators defined on the set of Boolean vectors $\{0,1\}^n$, given respectively by
\begin{align*}
  [F^+(x)]_i & := \min_{a \in A_i} \max_{b \in B_{i,a}} \max_{j:P_{i j}^{a b} > 0} x_j,\\
  [F^-(x)]_i & := \min_{a \in A_i} \max_{b \in B_{i,a}} \min_{j:P_{i j}^{a b} > 0} x_j,
\end{align*}
for every $i \in S$.
Note that these operators depend only on the support of the transition probabilities.

They suffice to characterize the families $\F^-$ and $\F^+$, as shown by the following lemma.
\begin{lemma}[\cite{AGH14}]
  \label{lem:BooleanReformulation}
  Let $F$ be a payment-free Shapley operator and let $I,J \subset S$.
  Thenn
  \begin{align*}
    I \in \F^- \enspace &\Leftrightarrow \enspace F^+(\indic_\Cpt{I}) \leq \indic_\Cpt{I},\\
    J \in \F^+ \enspace &\Leftrightarrow \enspace F^-(\indic_J) \geq \indic_J.
  \end{align*}
\end{lemma}

\begin{corollary}
  \label{coro:StructuralGaloisConnection}
  Given a payment-free Shapley operator, the families $\F^-$ and $\F^+$ and the Galois connection $(\Phi, \Phi^\star)$ are uniquely determined by the supports of the transition probabilities.
\end{corollary}

\section{Generic uniqueness of the bias vector of zero-sum games}
\label{sec:GenericUniqueness}

\subsection{The setting}
\label{sec-setting}

We fix the state space, the action spaces and the transition probabilities and consider parametric Shapley operators~\eqref{eq:ShapleyOperator}, that is, operators whose transition payments (thought of as parameters) can vary.
In particular, we denote by $\Rewards$ the {\em space of transition payments}, whose elements are the finite collections of reals $\big( r_i^{a b} \big)_{(i,a,b)}$ for the $(i,a,b) \in S\times A\times B$ such that $a\in A_i$ and $b\in B_{i,a}$.
This set can be identified with the Euclidean space $\R^q$ for a certain integer $q$.
For a given $r \in \Rewards$, we denote by $T_r$ the corresponding Shapley operator as defined in~\eqref{eq:ShapleyOperator}.

In this section we assume that the eigenproblem~\eqref{eq:NonlinearSpectralProblem} is solvable for all the values $r \in \Rewards$ of the transition payments.
For any Shapley operator $T$ as defined in~\eqref{eq:ShapleyOperator}, we denote by $\E(T)$ the set of its eigenvectors, and we say that $T$ has a unique eigenvector (up to an additive constant) if $\E(T)$ is reduced to a line (which is invariant by the addition of a constant).

Before stating the main result of this section, we recall that a polyhedral fan in $\R^q$ is a polyhedral complex consisting of polyhedral cones, that is, a finite set $\mathcal{K}$ of polyhedral cones such that
\begin{itemize}
\item $P \in \mathcal K$ and $F$ is a face of $P$ implies that $F \in \mathcal K$,
\item for all $P,Q \in \mathcal K$, $P \cap Q$ is a face of $P$ and $Q$.
\end{itemize}
We refer to the textbook~\cite{DLRS10} for background on polyhedral complexes and fans (but be aware that in this book a polyhedral complex is not necessarily finite).

\begin{theorem}
  \label{thm:GenericUniqueness}
  Assume that for any transition payment $r \in \Rewards$, the Shapley operator $T_r$ has an eigenvalue.
  Then the space $\Rewards$ is covered by a polyhedral fan such that for each $r$ in the interior of a full-dimensional cone, $T_r$ has a unique eigenvector.
\end{theorem}

Note that Theorem~\ref{thm:GenericUniqueness} implies that the transition payments $r \in \Rewards$ for which $T_r$ has more that one eigenvector (up to an additive constant) are situated in the finite union of subspaces of codimension $1$.

We prove Theorem~\ref{thm:GenericUniqueness} in the framework of both deterministic and stochastic games.
The structure of the proofs are similar but they rely on different theories, which we find interesting to present, although deterministic games can be seen as a particular case of stochastic games.
These theories are max-plus algebra and nonlinear Perron-Frobenius theory, respectively.

However, in both cases we infer from the following result a fundamental property about the geometry of the set of eigenvectors.
\begin{theorem}[\cite{Bru73}]
  Let $C$ be a closed convex subset of a finite dimensional Banach space.
  Let $T : C \to C$ be a nonexpansive map with respect to the norm of the Banach space.
  Then the fixed point set of $T$ is a nonexpansive retract of $C$, meaning that it is the image of a nonexpansive projection of $C$.
\end{theorem}

\begin{corollary}
  \label{coro:NonexpansiveRetract}
  Let $T: \R^n \to \R^n$ be an monotone and additively homogeneous map (hence nonexpansive in the sup-norm).
  Suppose Equation~\eqref{eq:NonlinearSpectralProblem} is solvable.
  Then $\E(T)$ is a retract of $\R^n$ by a nonexpansive map.
  In particular, $\E(T)$ is arcwise connected.
\end{corollary}

\subsection{Deterministic case}

In the deterministic case, each transition probability $P_i^{a b}$ has only one state in its support.
Thus, for $r \in \Rewards$, the operator $T_r$ in~\eqref{eq:ShapleyOperator} can be written as
\begin{equation}
  \label{eq:MinmaxFunction}
  [T_r(x)]_i = \min_{a \in A_i} \max_{1 \leq j \leq n} \big(\tilde r_{i j}^a + x_j \big),
\end{equation}
where $\tilde r_{i j}^a = \max \{r_i^{a b} \mid b \in B_{i,a}, \enspace P_{ij}^{a b}=1\}$ with the convention that $\max \emptyset = -\infty$.
Moreover, the set $C_{i,a}:=\set{j\in S}{\tilde r_{i j}^a\in \R}$ depends only on the transition probabilities, hence is independent of the value of the payments.

The above representation~\eqref{eq:MinmaxFunction} can be simplified using max-plus notations.
For this purpose, we introduce the max-plus semiring $\Rmax := \R \cup \{-\infty\}$, which is a commutative idempotent semiring with addition $x \oplus y = \max(x,y)$ and multiplication $x \otimes y = x+y$.
The zero and unit elements of $\Rmax$ are $-\infty$ and $0$, respectively.

Let $\Sigma$ be the (finite) set of policies of player {\sf MIN}, that is, the maps $\sigma: S \to A, i \mapsto \sigma(i) \in A_i$.
If $\R^n$ is endowed with its usual partial order then, for $r \in \Rewards$, the operator $T_r$ can be written as
\begin{equation}
  \label{eq:MaxPlusRepresentation}
  T_r(x) = \min_{\sigma \in \Sigma} M_r^\sigma x,
\end{equation}
where $M_r^\sigma \in \Rmax^{n \times n}$ is a max-plus matrix (that is a $n \times n$ matrix over the max-plus semiring $\Rmax$) and $M_r^\sigma x$ is to understand in the max-plus sense.
Precisely, we have $[M_r^\sigma]_{i j} = \tilde r^{\sigma(i)}_{i j}$.
Note that~\eqref{eq:MaxPlusRepresentation} means that, for every $x \in \R^n$, we have $T_r(x) \leq M_r^\sigma x$ for all $\sigma \in \Sigma$, and there is $\sigma \in \Sigma$ such that $T_r(x) = M_r^\sigma x$.
In particular, if $\lambda$ is the eigenvalue of $T_r$ and $u$ is an eigenvector, then there is a policy $\sigma$ such that $M_r^\sigma u = \lambda \unit + u$.

\subsubsection{Max-plus algebra}
We now present some results on max-plus matrices that can be found in~\cite[Chapter 3]{BCOQ92} or~\cite{ABG07}.
We draw the attention to the fact that vectors in $\Rmax^n$ are allowed to have entries equal to $-\infty$.

Let $M \in \Rmax^{n \times n}$ with no row identically equal to $-\infty$.
We need the following definitions.

\begin{definition}
  The {\em precedence graph} of $M$, denoted by $\G(M)$, is the weighted directed graph with set of nodes $S$ and an arc from $i$ to $j$ if $M_{i j} \neq -\infty$, in which case the weight of the arc is the real number $M_{i j}$.
\end{definition}

\begin{definition}
  The {\em maximal circuit mean}, $\rho(M)$, is the maximum average weight of the elementary circuits of $\G(M)$ (the circuits where every nodes are distinct, except the starting and the ending one).
\end{definition}

\begin{theorem}[{\cite[Rem. 3.24]{BCOQ92}}]
  Let $M \in \Rmax^{n \times n}$ with no row identically equal to $-\infty$.
  Then, there exists $u \in \Rmax^n$, with at least one finite entry, such that $ M u = \rho(M) \unit + u$.
\end{theorem}

The vocabulary is the same as for the eigenproblem~\eqref{eq:NonlinearSpectralProblem}: we say that $u$ is a tropical eigenvector of $M$ associated with the tropical eigenvalue $\rho(M)$.

\begin{definition}
  The {\em critical graph} of $M$, denoted by $\G^c(M)$, is the subgraph of $\G(M)$ consisting of all the elementary circuits with average weight equal to $\rho(M)$.
  A critical class of $M$ is a connected component of $\G^c(M)$.
\end{definition}

\begin{theorem}[{\cite[Thm. 3.101]{BCOQ92}}]
  \label{thm:TropicalEigenspace}
  Let $M \in \Rmax^{n \times n}$ with no row identically equal to $-\infty$.
  Denote by $\Emax(M)$ the set of tropical eigenvectors of $M$ associated with the eigenvalue $\rho(M)$.
  Then, $\Emax(M)$ is tropically generated by a finite family of vectors in $\Rmax^n$.
  Moreover, the least cardinality of a generating family of $\Emax(M)$ is equal to the number of critical classes of $M$.
\end{theorem}

\begin{theorem}
  Let $M \in \Rmax^{n \times n}$ and $\lambda \in \R$.
  If there is a vector $u \in \R^n$ such that $M u = \lambda \unit + u$, then $\lambda = \rho(M)$.
\end{theorem}

\subsubsection{Sketch of proof of Theorem~\ref{thm:GenericUniqueness}}
By definition of $\Rewards$, for every $\sigma$, the matrices $M_r^\sigma$ have the same precedence graph for all transition payments $r \in \Rewards$: there is an arc from $i$ to $j$ if $j \in C_{i,\sigma(i)}$.
Then, for every $\sigma$, we can define the piecewise linear function that maps an element $r \in \Rewards$ to the maximal circuit mean of $M_r^\sigma$, $\rho(M_rr^\sigma)$, and consider the polyhedral fan $\C^\sigma$ associated with the cone of linearity of this function.
To each cone of linearity with full dimension corresponds a unique critical circuit in $\G(M_r^\sigma)$ (where $r$ is any transition payment).

We consider then the polyhedral fan $\C$ obtained as the intersection of all the fans $\C^\sigma$.
We show that for each element $r \in \Rewards$ in the interior of a full-dimensional cone of $\C$, an eigenvector of $T_r$ is also a tropical eigenvector of a certain $M_r^\sigma$, associated with the eigenvalue $\rho(M_r^\sigma)$, and whose tropical eigenvector is unique (up to an additive constant).
This shows that  the eigenspace $\E(T_r)$ contains a finite number of elements up to the addition of constants. 
The conclusion follows from Theorem~\ref{coro:NonexpansiveRetract}.

\subsection{Stochastic case}
\label{sec-stochastic}

We go back to the general framework presented in Section~\ref{sec-setting}.
Recall that $\Sigma$ is the set of policies of player \MIN. 
For all $\sigma \in \Sigma$, we define $\Pi^\sigma$ as the set of policies of player \MAX, when the policy of \MIN\ is fixed to $\sigma$, that is the maps $\pi:S \to B, i \mapsto \pi(i) \in B_{i,\sigma(i)}$.
Moreover, given $\sigma \in \Sigma$ and $r \in \Rewards$,  we define the operator $T_r^\sigma:\R^n \to \R^n$ by
\[
  [T_r^\sigma(x)]_i = \max_{b \in B_{i,\sigma(i)}} \big( r_i^{\sigma(i) b} + P_i^{\sigma(i) b} x \big),
\]
for all $i\in S$.
Then, we have ($\R^n$ being endowed with its usual partial order)
\begin{equation}
  \label{eq:Tmin}
  T_r(x) = \min_{\sigma \in \Sigma} T_r^\sigma(x).
\end{equation}
Likewise, given $\sigma \in \Sigma$, we have
\begin{equation}
  \label{eq:StochasticControlOperator}
  T_r^\sigma(x) =  \max_{\pi \in \Pi^\sigma} \big( r^{\sigma \pi} + P^{\sigma \pi} x \big),
\end{equation}
where $r^{\sigma \pi}$ is the vector in $\R^n$ whose $i$-th entry is defined by $r^{\sigma \pi}_i = r^{\sigma(i) \pi(i)}_i$ and $P^{\sigma \pi}$ is the $n \times n$ stochastic matrix whose $i$-th row is given by $P^{\sigma \pi}_i = P^{\sigma(i) \pi(i)}_i$.

To prove Theorem~\ref{thm:GenericUniqueness}, we need to state some results on operators of type~\eqref{eq:StochasticControlOperator}, which are piecewise affine, convex (componentwise), monotone and additively homogeneous.

\subsubsection{Ergodic stochastic control}
Firstly, we characterize the eigenvalue of a componentwise convex operator $T^\sigma$ as defined in~\eqref{eq:StochasticControlOperator}.
\begin{lemma}
  \label{lem:EigenvalueStochasticControlOperator}
  Let $T^\sigma$ be an operator of type~\eqref{eq:StochasticControlOperator}.
  For every $\pi \in \Pi^\sigma$, let $m^{\sigma \pi}_C$ be the invariant measure associated with the final class $C$ of the stochastic matrix $P^{\sigma \pi}$.
  Denote by $\F^{\sigma \pi}$ the set of final classes of $P^{\sigma \pi}$.
  Suppose that $T^\sigma$ has an eigenvalue and denote it by $\lambda$.
  Then,
  \[
    \lambda = \max_{\pi \in \Pi^\sigma,\; C \in \F^{\sigma \pi}} \scalar{m^{\sigma \pi}_C}{r^{\sigma \pi}}.
  \]
\end{lemma}

Secondly, we characterize the set of eigenvectors of the operator $T^\sigma$.
For that, we need few definitions that are detailed in~\cite{AG03}.
For a map $F:\R^n \to \R^n$, componentwise convex, the notion of subdifferential is generalized by setting, for $x \in \R^n$, $\partial F(x) = \{P \in \R^{n \times n} \mid P(y-x) \leq F(y)-F(x), \enspace \forall y \in \R^n\}$.
If $F$ is monotone and additively homogeneous, it can be shown that $\partial F(x)$ is a convex set of stochastic matrices.
Suppose that $F$ has an eigenvector $u$.
Then the critical graph of $F$, denoted by $\G^c(F)$, is the directed graph defined as the union of the final graphs of the stochastic matrices $P \in \partial F(u)$.
This graph does not depend on the choice of the eigenvector $u$.
A critical node of $F$ is a node of $\G^c(F)$ and a critical class of $F$ is the set of nodes of a strongly connected component of $\G^c(F)$.
\begin{theorem}[{\cite[Thm 1.1]{AG03}}]
  \label{thm:ConvexSpectralTheorem}
  Let $F:\R^n \to \R^n$ be a convex, monotone and additively homogeneous map.
  Suppose that $F$ has an eigenvector and denote by $C$ the set of its critical nodes.
  Then,
  (i) the restriction $\mathsf{r}:\R^n \to \R^C, x \mapsto (x_i)_{i \in C}$ is an isomorphism from $\E(F)$ to its image $\E^c(F)$;
  (ii) $\E^c(F)$ is a convex set whose dimension is at most equal to the number of critical classes of $F$, and this bound is attained when $F$ is piecewise affine.
\end{theorem}

\subsubsection{Sketch of proof of Theorem~\ref{thm:GenericUniqueness}}
Let $r \in \Rewards$.
For every $\sigma$, the operaor $T_r^\sigma$ is convex, monotone and additively homogeneous.
Then, for every $\sigma$, we can define the piecewise linear map $\rho^\sigma: r \in \Rewards \mapsto \max_{\pi \in \Pi^\sigma,\; C \in \F^{\sigma \pi}} \scalar{m^{\sigma \pi}_C}{r^{\sigma \pi}}$, where, for $\pi \in \Pi^\sigma$, $m^{\sigma \pi}_C$ is the invariant measure associated with the final class $C$ of the stochastic matrix $P^{\sigma \pi}$.
We consider the polyhedral fan $\C^\sigma$ associated with the cone of linearity of the map $\rho^\sigma$.
To each cone of linearity with full dimension corresponds a unique final class $C$.

We consider then the polyhedral fan $\C$ obtained as the smallest refinement of all the fans $\C^\sigma$.
We show that for each transition payment $r$ in the interior of a full-dimensional cone of $\C$, an eigenvector $u$ of $T_r$ is also an eigenvector of a certain $T_r^\sigma$ with a unique critical class, hence a unique eigenvector (up to an additive constant).
This shows that the eigenspace $\E(T_r)$  contains a finite number of elements up to the addition of constants. 
The conclusion follows from Theorem~\ref{coro:NonexpansiveRetract}.

\section{Application to policy iteration}
\label{sec:Application}

Let us use the notations of Section~\ref{sec-stochastic}.
In particular, $\Sigma$ is the (finite) set of policies of player \MIN, that is the maps $\sigma: S \to A, i \mapsto \sigma(i)\in A_i$, $\Pi^\sigma$  is the (finite) set of policies of player \MAX\ when the policy of \MIN, $\sigma \in \Sigma$, is fixed, that is the set of maps $\pi:S \to B, i \mapsto \pi(i) \in B_{i,\sigma(i)}$, and for all $\sigma\in\Sigma,\; \pi\in \Pi^\sigma$, $P^{\sigma \pi}$ is the $n \times n$ stochastic matrix whose $i$th row is given by $P^{\sigma \pi}_i = P^{\sigma(i) \pi(i)}_i$.

When $T$ is a Shapley operator~\eqref{eq:ShapleyOperator} with finite state and action spaces (Assumption~\ref{StandardAssumptions}) such that  all the stochastic matrices $P^{\sigma\tau}$ are irreducible, Hoffman and Karp~\cite{HK66} introduced the following algorithm with input $T$ and output an eigenvalue $\lambda$ and a bias vector $u$ of $T$.
Then optimal stationary policies for both players can be computed.

\begin{algo}[\cite{HK66}]
  \label{algo-main}
  Select an arbitrary policy $\sigma_0\in\Sigma$ of player \MIN, then apply successively  the two following steps for $\ki\geq 0$, until  $\sigma_{\ki+1}=\sigma_{\ki}$ :
  \begin{enumerate}
    \item \label{valuee} Compute the eigenvalue $\lambda^{\ki}$ and the bias  $v^{\ki}$ of the game with fixed policy $\sigma_\ki$ of player \MIN, that is the solutions $\lambda$ and  $v$ of $\lambda \unit +v = T^{\sigma_{\ki}}(v)$;
    \item Improve the policy: choose an optimal policy for  $v^{\ki}$, that is $\sigma_{\ki+1}\in\Sigma$  such that $T(v^{\ki})=T^{\sigma_{\ki+1}}(v^{\ki})$, with $\sigma_{\ki+1}=\sigma_{\ki}$ as soon as this is possible.
  \end{enumerate}
  Return $\lambda^{\ki}$ and $v^{\ki}$.
\end{algo}
In the above algorithm, Step~\ref{valuee}) is solved by using the same policy iteration for the (one-player) game with fixed policy $\sigma_{\ki}$, which constructs $\lambda^{\ki,\kj}$, $v^{\ki,\kj}$ and $\pi_{\ki,\kj}$ from some $\pi_{\ki,0}$.
One can show that the sequence $\lambda^{\ki}$ of Algorithm~\ref{algo-main} is nonincreasing ($\lambda^{\ki+1}\leq \lambda^{\ki}$), and that when $\lambda^{\ki+1}= \lambda^{\ki}$, then $v^{\ki+1}=v^{\ki}$ (up to an additive constant) and $\sigma_{\ki+1}=\sigma_{\ki}$.
Hence, the algorithm terminates after a finite number of steps, for the set of policies $\Sigma$ is finite.

When the stochastic matrices $P^{\sigma\tau}$ are not irreducible, the existence of $\lambda^{\ki}$ and $v^{\ki}$ may fails in Step~\ref{valuee}) of Algorithm~\ref{algo-main}.
Moreover, even if $\lambda^{\ki}$ and $v^{\ki}$ can be constructed at each step $\ki$, the nonuniqueness of $v^{\ki}$ (up to an additive constant) may lead to a cycling of the sequence $v^{\ki}$ in Algorithm~\ref{algo-main}, as shown for instance in~\cite[Section 6]{ACTDG12}.
In the general case, one need first to replace eigenvectors by invariant half-lines to avoid problems of existence and second to apply a special treatment of degeneracies (when the eigenvalue or mean payoff vector does not change from one step to another and the bias vector or invariant half-line is not unique) as in~\cite{cochet-cras,ACTDG12,bourque} to avoid cycling of the algorithm.

Using the arguments of Section~\ref{sec-stochastic}, we obtain the following result, under the same setting as for Theorem~\ref{thm:GenericUniqueness}, that is fixing the state space, the action spaces and the transition probabilities and considering parametric Shapley operators~\eqref{eq:ShapleyOperator}.
\begin{theorem}
  \label{thm:GenericPolicy}
  The space $\Rewards$ is covered by a polyhedral fan such that for each transition payment $r$ in the interior of a full-dimensional cone, if Algorithm~\ref{algo-main} is well posed, that is if at each step $\ki$ of the algorithm, the operator $T_r^{\sigma_{\ki}}$ has an eigenvalue, then it terminates after a finite number of steps and gives an eigenvalue $\lambda$ and a bias vector $u$ of $T_r$.
\end{theorem}

{\noindent\hspace{1em}{\itshape Sketch of proof:}}
Considering the same fan $\C$ as in Section~\ref{sec-stochastic}, we get that for each $r$ in the interior of a full-dimensional cone of $\C$, and for each $\sigma \in \Sigma$, if the operator $T_r^\sigma$  has an eigenvector $u$, then this eigenvector is unique (up to an additive constant).
This implies that, if at each step $\ki$ of the algorithm, the operator $T_r^{\sigma_{\ki}}$ has an eigenvalue $\lambda^{\ki}$, then the bias vector $v^{\ki}$ is unique up to an additive constant.
One can always show (without any assumption) that the sequence $\lambda^{\ki}$  is nonincreasing ($\lambda^{\ki+1}\leq \lambda^{\ki}$).
Then, using the uniqueness of bias vectors, we deduce from Lemma 3.3 of~\cite{AG03}, that if $\lambda^{\ki+1}= \lambda^{\ki}$, then up to an additive constant, $v^{\ki+1} \leq v^{\ki}$ with equality on the set of critical nodes of $T_r^{\sigma_{\ki}}$.
This implies that the sequence $v^{\ki}$ coincides up to an additive constant with the sequence obtained in the algorithm introduced in~\cite{cochet-cras} and developed in~\cite{ACTDG12}.
Then, the convergence of Algorithm~\ref{algo-main} in a finite number of steps follows from the convergence of the algorithm of~\cite{cochet-cras,ACTDG12}.

\section{Example}

Consider the following parametric Shapley operator defined on $\R^3$ (here we use $\wedge$ and $\vee$ instead of $\min$ and $\max$, respectively, and we recall that $+$ has precedence over them):
\begin{equation*}
\small
T_r(x) =
\begin{pmatrix}
{\color{red} r_1} + \frac{1}{2} (x_1 + x_3) \, \wedge \, {\color{red} r_2} + \frac{1}{2} (x_1 + x_2)\\
{\color{red} r_3} + \frac{1}{2} (x_1 + x_3) \, \wedge \, \left( {\color{red} r_4} + \frac{1}{2} (x_1 + x_2) \, \vee \, {\color{red} r_5} + x_3 \right)\\
{\color{red} r_6} + \frac{1}{2} (x_1 + x_3) \, \vee \, {\color{red} r_7} + x_3
\end{pmatrix}
\end{equation*}
with $r \in \R^7$.
Its recession function is
\begin{equation*}
\small
\hat T_r(x) =
\begin{pmatrix}
\frac{1}{2} (x_1 + x_3) \, \wedge \, \frac{1}{2} (x_1 + x_2)\\
\frac{1}{2} (x_1 + x_3) \, \wedge \, \left( \frac{1}{2} (x_1 + x_2) \, \vee \, x_3 \right)\\
\frac{1}{2} (x_1 + x_3) \, \vee \, x_3
\end{pmatrix}
\end{equation*}
and it can be checked that it has only trivial fixed points.
Hence, according to Theorem~\ref{thm:StructuralEigenvalue}, Equation~\eqref{eq:NonlinearSpectralProblem} is solvable for every $r \in \R^7$.
Alternatively, we may consider the Galois connection between the families
$\F^- = \big\{ \emptyset, \{1,3\}, \{1,2,3\} \big\}$ and $\F^+ = \big\{ \emptyset, \{3\}, \{1,2,3\} \big\}$ as defined in~\ref{sec:GaloisConnection}.
We  verify that the only nontrivial element of $\F^-$ is not closed with respect to this Galois connection (since $\Phi^\star \circ \Phi (\{1,3\}) = \Phi^\star(\emptyset) = \{1,2,3\}$), which leads to the same conclusion according to Corollary~\ref{coro:GaloisConnectionFixedPoint}.

We now fix $r=r_0=(0,1,2,1,-2,-3,1)$ and denote by $T_0$ the corresponding operator.
We consider the perturbations of $r_0$ that depend only on the state i.e., in the space of payments $\R^7$, we consider the affine subspace $\mathcal{A}_0 = \{r_0+(g_1,g_1,g_2,g_2,g_2,g_3,g_3) \mid g \in \R^3 \}$, corresponding to the operators $T_0 + g$ with $g \in \R^3$ (note that by additive homogeneity of $T_r$, we may assume w.l.o.g.~ that $g_3=0$).
Figure~\ref{fig} shows the polyhedral complex obtained as the intersection between the polyhedral fan of Theorem~\ref{thm:GenericUniqueness} and the affine subspace $\mathcal A_0$.
Here for each $g$ in the interior of a full-dimensional polyhedron, $T_0+g$ has a unique eigenvector (up to an additive constant).
Let us detail what happens in the neighborhood of $g=0$, point in which $T_0+g$ fails to have a unique eigenvector.
Note that in the neighborhood of $g=0$, the eigenvalue of $T_0+g$ remains $1$.

\begin{figure}[!h]

\begin{minipage}{0.4\linewidth}
\begin{center}
\begin{tikzpicture}[scale=0.13]
\draw [->, >= angle 60,gray,thick] (-5,0) -- (20,0);0
\draw (19,0) node[above] {$g_1$};
\foreach \x in {0,5,10,15}
	\draw [gray,thick] (\x,0) -- (\x,-0.5);
\draw (10,0) node[below] {\footnotesize $10$};
\draw [->, >= angle 60,gray,thick] (0,-20) -- (0,7);
\draw (0,6) node[right] {$g_2$};
\foreach \y in {-15,-10,-5,0,5}
	\draw [gray,thick] (-0.5,\y) -- (0,\y);
\draw (0,-10) node[left] {\footnotesize $-10$};
\draw (15,6) node {$g_3=0$};
\draw [very thick,orange] (5,-20) -- (5,7);
\draw [very thick,blue] (-2,7) -- (23/2,-20);
\draw [very thick,red] (11,-11) -- (-5,5);
\draw [very thick,red] (11,-11) -- (20,-42/3);
\draw [very thick,red] (11,-11) -- (16,-20);
\draw (0,0) node {$\bullet$};
\draw (0,0) node [below left] {$0$};
\end{tikzpicture}
\end{center}
\caption{}
\label{fig}
\end{minipage}
\hfill
\begin{minipage}{0.55\linewidth}
If $g_1+g_2 = 0$, the eigenvectors of $T_0+g$ are defined by $x_1=x_2+2 g_1$ and $-3+g_2 \leq x_2-x_3 \leq -2-g_1$.

If $g_1+g_2 > 0$, the unique eigenvector (up to an additive constant) is $(-2+2 g_1,-2+2 g_1+2 g_2,0)$.

If $g_1+g_2 < 0$, the unique eigenvector is $(-3+2 g_1+g_2,-3+g_2,0)$.
\end{minipage}

\vspace{-1em}
\end{figure}

\addtolength{\textheight}{-22cm}   






\bibliographystyle{IEEEtran}
\bibliography{references}

\end{document}